\magnification=1200

\input amstex

\documentstyle{amsppt}


\hsize=165truemm

\vsize=227truemm


\def\p#1{{{\Bbb P}^{#1}_{k}}}

\def\a#1{{{\Bbb A}^{#1}_{k}}}

\def\Hilb{{{\Cal H}\kern -0.25ex{\italic ilb\/}}}

\def\Scand{{{\Cal S}\kern -0.25ex{\italic cand\/}}}

\def\Hom{{{\Cal H}\kern -0.25ex{\italic om\/}}}

\def\Ext{{{\Cal E}\kern -0.25ex{\italic xt\/}}}

\def\Sim{{{\Cal S}\kern -0.25ex{\italic ym\/}}}

\def\Ker{{{\Cal K}\kern -0.25ex{\italic er\/}}}

\def\GL{\operatorname{GL}}

\def\PGL{\operatorname{PGL}}

\def\Sing{\operatorname{Sing}}

\def\emdim{\operatorname{emdim}}

\def\char{\operatorname{char}}

\def\gr{\operatorname{gr}}

\def\spec{\operatorname{spec}}

\def\proj{\operatorname{proj}}

\def\ext{\operatorname{Ext}}

\def\rk{\operatorname{rk}}

\def\lev{\operatorname{lev}}

\def\Soc{\operatorname{Soc}}

\def\Ofa#1{{{\Cal O}_{#1}}}

\def\M{\operatorname{\frak M}}

\def\ga#1{{\accent"12 #1}}


\topmatter

\title
On the Gorenstein locus of some punctual Hilbert schemes
\endtitle


\author
Gianfranco Casnati, Roberto Notari
\endauthor

\address
Gianfranco Casnati, Dipartimento di Matematica, Politecnico di Torino,
c.so Duca degli Abruzzi 24, 10129 Torino, Italy
\endaddress

\email
casnati\@calvino.polito.it
\endemail

\address
Roberto Notari, Dipartimento di Matematica \lq\lq Francesco Brioschi\rq\rq, Politecnico di Milano,
via Bonardi 9, 20133 Milano, Italy
\endaddress

\email
roberto.notari\@polimi.it
\endemail

\keywords
Hilbert scheme, arithmetically Gorenstein subscheme, Artinian algebra
\endkeywords

\subjclassyear{2000}
\subjclass
14C05, 13H10, 14M05
\endsubjclass

\abstract
Let $k$ be an algebraically closed field and let $\Hilb_{d}^{G}(\p{N})$ be the open locus of the Hilbert scheme $\Hilb_{d}(\p{N})$ corresponding to Gorenstein subschemes. We prove that $\Hilb_{d}^{G}(\p{N})$ is irreducible for $d\le9$, we characterize geometrically its singularities for $d\le 8$ and we give some results about them when $d=9$ which give some evidence to a conjecture on the nature of the singular points in $\Hilb_{d}^{G}(\p{N})$. 
\endabstract

\endtopmatter

\document

\head
1. Introduction and notation
\endhead

Let $k$ be an algebraically closed field and denote by $\Hilb_{p(t)}(\p N)$ the Hilbert scheme parametrizing closed subschemes in $\p N$ with fixed Hilbert polynomial $p(t)\in{\Bbb Q}[t]$. Since A\. Grothendieck's proof  in [Gr] of the existence of $\Hilb_{p(t)}(\p N)$, the problem of finding a useful description of this scheme has attracted the interest of many researchers in algebraic geometry.

One of the first, now well known, results in this direction is due to R\. Hartshorne who proved the connectedness of $\Hilb_{p(t)}(\p N)$ in [Ha1]. There have also been studies of some loci in the Hilbert scheme. Two of the most notable are the description of the locus of codimension $2$ arithmetically Cohen--Macaulay subschemes (see [El]) and of the locus of codimension $3$ arithmetically Gorenstein subschemes (see [MR] and [JK--MR]).

In the study of punctual Hilbert schemes, the first fundamental result is due to J\. Fogarty who proved the irreducibility and smoothness of  $\Hilb_{d}(\p N)$, $d\in\Bbb N$, when $N=2$. The same result holds more generally if one considers subschemes of codimension $2$ of any smooth surface (see [Fo]).

In [Ia1] the author proved that, if $d$ is large with respect to $N$, $\Hilb_{d}(\p N)$ is never irreducible. 
Indeed for every $d$ and $N$ there always exists a generically smooth component of $\Hilb_{d}(\p N)$ having dimension $dN$ whose general point corresponds to a reduced set of $d$ points but, for $d$ large with respect to $N>2$, there is at least one other component whose general point corresponds to an irreducible scheme of degree $d$ supported on a single point. 

In view of these earlier works it is reasonable to consider the irreducibility and smoothness of other naturally occuring loci in $\Hilb_{d}(\p N)$. E.g: one of the loci that has interested us is the set $\Hilb_{d}^G(\p N)$ of points in $\Hilb_{d}(\p N)$ representing schemes which are Gorenstein. This is an important locus since it includes reduced schemes. 

Some results about $\Hilb_{d}^G(\p N)$ are known. E.g. since $\Hilb_{d}^G(\p N)$ contains all reducible schemes of degree $d$ it follows that $\Hilb_{d}^G(\p N)$ contains an open (not necessarily dense) subset of $\Hilb_{d}(\p N)$. More precisely $\Hilb_{d}^G(\p N)$ is actually open since its complement in $\Hilb_{d}(\p N)$ coincides with the projection on $\Hilb_{d}(\p N)$ locus of point over which the relative dualizing sheaf of the universal flat family is not invertible. Another result, part of the folklore, gives the irreducibility and smoothness of $\Hilb_{d}^G(\p N)$ when $N=3$. We provide a proof of this fact in Section 5. In [I--E] and [I--K] it is shown that $\Hilb_{d}^G(\p N)$ is never irreducible for $d\ge 14$ and $N\ge6$.

These results leave open the question of irreducibility for small $d$ and all $N>3$. We are aware of some scattered results about the existence of singular points on $\Hilb_{d}(\p N)$ can be found in [Ka], [C--N]. One of our principal results on these matters is the following theorems proved in Section 5.

\proclaim{Theorem A}
Assume the characteristic of $k$ is $p\ne2,3$. The locus $\Hilb_{d}^{G}(\p{N})$ is irreducible for $d\le 9$. 
\qed
\endproclaim

Very recently, in [C--E--V--V], the authors prove both the irreducibility of  $\Hilb_{d}(\p{N})$ when $d\le7$ and the existence of exactly two components in  $\Hilb_{8}(\p{N})$, $N\ge4$.

In order to prove Theorem A we need to study deformations of some particular local Artinian Gorenstein $k$--algebras of degree $d\le 9$ and embedding dimension at least $4$. We begin the study of such algebras in Section 2 where we fix the notation and recall some elementary facts. In the final part of Section 2 and in Sections 3 and 4 we give a complete classification of such kind of algebras.

The problem of classifiying  local Artinian $k$--algebras is classical. It is completely solved for $d\le6$ (see [Ma1], [Ma2], [C--N] when $\char(k)>3$) and [Po2] without any restriction on the characteristic. When $d\ge7$ it is classically known that such algebras have moduli and their parameter spaces have been the object of deep study (again see [Ma2], [I--E] and [Po1]).

Let us now restrict to the {\sl Gorenstein} $k$--algebras. Their classification in degree $d=7$ follows from the result classically proved in [Sa] (see also [C--N]) and, when the characteristic of $k$ is $p=0$, in the more recent paper [E--V]. In degrees $d=8,9$ a complete classification can be done using, in addition to those papers, also the work  [Wa] and [E--I] on the classification of nets of conics.

Then we turn our attention to the singularities of the Hilbert scheme. In Section 5 we also prove the following

\proclaim{Theorem B}
If $d\le8$, then $X\in\Sing(\Hilb_{d}^{G}(\p{N}))$  if and only if the corresponding scheme $X$ has embedding dimension $4$ at least at one of its points. 
\qed
\endproclaim

This is no longer true if $d=9$. We examine, at the end of Section 5, this case and we give a partial result which gives some evidence to a conjeture on $\Sing(\Hilb_{d}^{G,gen}(\p{N}))$ for each $d$.

In order to prove the above theorem we have to combine on one hand the study the hierachy of local Artinian Gorenstein $k$--algebras of degree $d\le 9$ and embedding dimension at least $4$, on the other the properties of the $G$--fat point $X\subseteq\p{d-2}$ proved in [C--N].

We would like to express our thanks to A. Conca, A. Geramita, A. Iarrobino and J.O. Kleppe for some interesting and helpful suggestions. A particular thank goes to J. Elias and G. Valla, who pointed out an incongruence in the results proved in Section 3.1 with respect to the more general classification of almost stretched local, Artinian, Gorenstein $k$--algebras described  in [E--V].

\subhead
Notation
\endsubhead
In what follows $k$ is an algebraically closed field. We denote its characteristic by $\char(k)$. 

Recall that a Cohen--Macaulay local ring $R$ is one for which
$\dim(R)={\roman{depth}}(R)$. If, in addition the injective dimension of $R$ is finite then $R$ is called Gorenstein (equivalently, if $\ext_R^i\big(M,R)=0$ for each $R$--module $M$ and $i>\dim(R)$).  An
arbitrary ring $R$ is called Cohen--Macaulay (resp. Gorenstein) if
$R_{\frak M}$ is Cohen--Macaulay (resp. Gorenstein) for every
maximal ideal ${\frak M}\subseteq R$. 

All the schemes $X$ are separated and of finite type over $k$. A scheme $X$ is Cohen--Macaulay (resp. Gorenstein) if for each
point $x\in X$ the ring ${\Cal O}_{X,x}$ is
Cohen--Macaulay (resp. Gorenstein). The scheme $X$ is Gorenstein if and only if it is Cohen--Macaulay and its
dualizing sheaf $\omega_{X}$ is invertible.

For each numerical polynomial $p(t)\in{\Bbb Q}[t]$ of degree at most $n$ we denote by $\Hilb_{p(t)}(\p N)$ the Hilbert scheme of closed subschemes of $\p N$ with Hilbert polynomial $p(t)$. With abuse of notation we will denote by the same symbol both a point in $\Hilb_{p(t)}(\p N)$ and the corresponding subscheme of $\p N$. Moreover we denote by $\Hilb_{p(t)}^G(\p N)$ the locus of points representing Gorenstein schemes. This is an open subset of $\Hilb_{p(t)}(\p N)$, though not necessarily dense.

If $X\subseteq\p N$ we will denote by $\Im_X$ its sheaf of ideals in $\Ofa X$ and we define the normal sheaf of $X$ in $\p N$ as ${\Cal N}_X:=\Hom_{X}\big(\Im_{X}/\Im_{X}^2,\Ofa{X}\big)$. The homogeneous ideal of $X$ is
$$
I_X:=\bigoplus_{t\in\Bbb Z}H^0\big(X,\Im_X(t)\big)\subseteq\bigoplus_{t\in\Bbb Z}H^0\big(\p N,\Ofa{\p N}(t)\big)\cong S:=k[x_0,\dots,x_N].
$$
The ideal $I_X$ is saturated. We define the homogeneous coordinate ring as $S_X:=S/I_X$. We have $X=\proj(S_X)$ and the embedding $X\subseteq \p N$ corresponds to the canonical epimorphism $S\twoheadrightarrow S_X$. The scheme $X$ is said arithmetically Gorenstein (briefly aG) if $S_X$ is a Gorenstein ring.

\head
2.  The locus $\Hilb_{d}^{G}(\p{N})$
\endhead

As explained in the introduction we denote by  $\Hilb_{d}^{G}(\p{N})\subseteq \Hilb_{d}(\p{N})$ the Gorenstein locus, i.e. the locus of points in $\Hilb_{d}(\p{N})$ representing Gorenstein subschemes of $\p N$. The locus $\Hilb_{d}^G(\p{N})$ is open, but is not necessarily dense, in $\Hilb_{d}(\p{N})$.

Reduced schemes obviously represent points in $\Hilb_{d}^{G}(\p{N})$. Clearly the locus of such schemes in $\Hilb_{d}(\p{N})$ is birational to a suitable open subset of the $d$--th symmetric product of $\p N$, thus is irreducible of dimension $dN$ (see [Ia1]) and we will denote by $\Hilb_{d}^{gen}(\p{N})$ its closure in $\Hilb_{d}(\p{N})$. It follows that $\Hilb_{d}^{G,gen}(\p{N}):=\Hilb_{d}^G(\p{N})\cap\Hilb_{d}^{gen}(\p{N})$ is irreducible of dimension $dN$ and open in $\Hilb_{d}^G(\p{N})$.

Let $X\in \Hilb_{d}^G(\p{N})$. Then $X=\bigcup_{i=1}^p X_i$ where the $X_i$ are irreducible and pairwise disjoint of degree $d_i$, with $d=\sum_{i=1}^pd_i$.

Fix one of such component and call it $Y$. Each such scheme is affine, say $Y\cong\spec(A)$ where $A$ is an Artinian, Gorenstein $k$--algebra of degree $\delta$, i.e. with $\dim_k(A)=\delta$, and maximal ideal $\M$. In order to study our scheme $Y$, hence $X$, it is then natural to study $A$.

Let the embedding dimension of $A$, i.e. $\emdim(A):=\dim_k(\M/\M^2)$, be at most $n\le N$. Then we have a surjective morphism from the symmetric $k$--algebra on $\M/\M^2$, which is $k[y_1,\dots,y_n]$, onto $A$. Hence we have an isomorphism $A\cong k[y_1,\dots,y_n]/I$. Assume that $Y$ does not intersect the hyperplane $\{\ x_0=0\ \}$. Then the embedding $Y\subseteq\p{N}$ corresponds to an epimorphism $k[x_1,\dots,x_{N}]\twoheadrightarrow k[y_1,\dots,y_n]/I$. Due to the definition of $n$ such a morphism factors through another epimorphism $k[x_1,\dots,x_{N}]\twoheadrightarrow k[y_1,\dots,y_n]$, defining a subscheme of $\a{N}$ isomorphic to $\a n$. Its closure $Q$ in $\p{N}$ is obviously smooth around $Y$. 

Thus each flat family ${\Cal Y}\subseteq Q\times B$ with special fibre $Y$ lifts to a flat family ${\Cal Y}\subseteq \p N\times B$. We conclude that if each component of $X$ is smoothable in $\a n$, then $X\in\Hilb_{d}^{G,gen}(\p{N})$.

Taking into account Corollary 4.3 of [HK], this holds in particular for schemes which are AS in the sense of the following

\definition{Definition 2.1}
Let $X$ be a scheme of dimension $0$. We say that $X$ is AS (almost solid) if the embedding dimension at every point of $X$ is at most three.
\enddefinition

Now we turn our attention to the singular locus of $\Hilb_{d}^{G,gen}(\p{N})$. Let $X=\bigcup_{i=1}^p X_i\in \Hilb_{d}^{G,gen}(\p{N})$ be as above. Since
$h^0\big(X,{\Cal N}_{X}\big)=\bigoplus_{i=1}^ph^0\big(X_i,{\Cal N}_{X_i}\big)$ and $h^0\big(X_i,{\Cal N}_{X_i}\big)\ge d_iN$, it turns out that $X$ is obstructed if and only if the same is true for at least one of the $X_i$.

Again, from now on, we will fix our attention on the above irreducible $Y\cong\spec(A)\in \Hilb_{\delta}^{G}(\p{N})$. Consider the standard sequence
$$
0\longrightarrow\Im_Q\longrightarrow\Im_Y\longrightarrow\Im_{Y\vert Q}\longrightarrow0,
$$
where $\Im_{Y\vert Q}\subseteq\Ofa Q$ is the sheaf of ideals of $Y\subseteq Q$. Since $\Hom_{\p{N}}\big({\Cal F},\Ofa{Y}\big)\cong\Hom_{Y}\big({\Cal F}\otimes\Ofa Y,\Ofa{Y}\big)$ for each sheaf $\Cal F$ of $\Ofa Y$--modules, then we have the exact sequence
$$
0\longrightarrow{\Cal N}_{Y\vert Q}\longrightarrow{\Cal N}_{Y}\longrightarrow\Hom_{Y}\big(\Im_Q\otimes\Ofa Y,\Ofa{Y}\big)\longrightarrow\Ext_{Y}^1\big(\Im_{Y\vert Q}/\Im_{Y\vert Q}^2,\Ofa{Y}\big),
$$
where ${\Cal N}_{Y\vert Q}:=\Hom_{Y}\big(\Im_{Y\vert Q}/\Im_{Y\vert Q}^2,\Ofa{Y}\big)$. Since $Q$ is locally complete intersection around $Y$ we have, in a suitable open neighbourhood ${\Cal U}\subseteq\p N$ of $Y$, an exact Koszul complex of the form
$$
\wedge^2{\Cal O}_{\Cal U}^{\oplus N-n}\longrightarrow{\Cal O}_{\Cal U}^{\oplus N-n}\longrightarrow{\Im_Q}_{\vert\Cal U}\longrightarrow 0.
$$
whence $\Im_Q\otimes\Ofa Y\cong{\Cal O}_Y^{\oplus N-n}$.
It follows the existence of an exact sequence
$$
0\longrightarrow{\Cal N}_{Y\vert Q}\longrightarrow{\Cal N}_{Y}\longrightarrow{\Cal O}_Y^{\oplus N-n}\longrightarrow\Ext_{Y}^1\big(\Im_{Y\vert Q}/\Im_{Y\vert Q}^2,\Ofa{Y}\big)
$$
Since $Y$ is Gorenstein $\Ext_{Y}^1\big(\Im_{Y\vert Q}/\Im_{Y\vert Q}^2,\Ofa{Y}\big)=0$. Thus, taking cohomology, we obtain
$$
h^0\big(Y,{\Cal N}_{Y}\big)= h^0\big(Y,{\Cal N}_{Y\vert Q}\big)+(N-n)h^0\big(Y,\Ofa Y\big)=h^0\big(Y,{\Cal N}_{Y\vert \a n}\big)+(N-n)\delta.\tag 2.2
$$
Since $h^0\big(Y,{\Cal N}_{Y\vert \a n}\big)\ge n\delta$, it follows that $h^0\big(Y,{\Cal N}_{Y}\big)\ge N\delta$.

In particular the obstructedness of a smoothable $X\in \Hilb_{d}^{G}(\p{N})$ does not depend on its embedding but only on the obstructedness of its irreducible components in the space of lower dimension which contains them. Taking into account Proposition 2.2 and Remark 2.3 of [JK--MR], we then deduce the unobstructedness of AS Gorenstein schemes of dimension $0$.

We can summarize the above discussion in the following

\proclaim{Proposition 2.3}
Let $\char(k)\ne2$. If $X\in\Hilb_d^G(\p N)$ represents an AS scheme, then $X\in\Hilb_d^{G,gen}(\p N)$ and it is unobstructed.
\qed
\endproclaim

For the reader's benefit we recall the following well--known

\proclaim{Corollary 2.4}
Let $\char(k)\ne2$. If $N\le3$ then $\Hilb_d^G(\p N)$ is irreducible and smooth.
\endproclaim
\demo{Proof}
With the above hypotheses we have $\Hilb_d^G(\p N)=\Hilb_d^{G,gen}(\p N)$, which is then irreducible.
\qed
\enddemo

It is then natural to ask if $\Hilb_{d}^{G}(\p{N})$ is irreducible. Or, equivalently, are non--AS schemes smoothable, i.e. in $\Hilb_{d}^{G,gen}(\p{N})$?

The answer to this question is in general negative. As pointed out in the introduction, in Section 6.2 of [I--K] the authors states the existence of local Artinian, Gorenstein $k$--algebras $A$ the deformations of which are all of the same type, using a method previously introduced in [I--E]: thus such kind of algebras $A$ define an irreducible component in $\Hilb_{14}^{G}(\p{6})$ distinct from $\Hilb_{14}^{G,gen}(\p{6})$.

A second natural question is to ask if $\Sing(\Hilb_{d}^{G}(\p{N}))$ coincides with the locus of non--AS schemes

When $d$ increases it is the answer to the above question is again negative. E.g. take $X:=\spec(A)\in\Hilb_{16}^{G}(\p{4})$, where $A:=k[x_1,x_2,x_3,x_4]/(x_1^2,x_2^2,x_3^2,x_4^2)$. Thus $X$ is a complete intersection, thus it is trivially smoothable, hence it belongs to the component $\Hilb_{16}^{G,gen}(\p{4})$ which has dimension $64$. Being a complete intersection, ${\Cal N}_X$ is locally free, thus it is actually free, since $X$ has dimension $0$. This means that ${\Cal N}_X\cong{\Cal O}_X^{\oplus4}$, thus $h^0\big(X{\Cal N}_X)=4h^0\big(X,{\Cal O}_X)=64$, hence $X$ turns out to be unobstructed.

The object of our paper is to prove the irreducibility of $\Hilb_{d}^{G}(\p{N})$ and to characterize its singular locus , when $d\le9$. Due to Theorem 2.3 above it is clear that we have to focus our attention on non--AS schemes $X\in \Hilb_{d}^G(\p{N})$.
To this purpose we first look at the intrinsic structure of local Artinian, Gorenstein $k$--algebras of degree $d\le9$ and $\emdim(A)\ge4$. 

Let $A$ be a local Artinian $k$--algebra of degree $d$ with maximal
ideal $\M$. 
In general we have a filtration
$$
A\supset\M\supset\M^2\supset\dots\supset\M^e\supset\M^{e+1}=0 
$$
for some integer $e\ge1$, so that its associated graded algebra
$$
\gr(A):=\bigoplus_{i=0}^\infty\M^i/\M^{i+1}
$$
is a vector space over $k\cong A/\M$ of finite dimension $d=\dim_k(A)=\dim_k(\gr(A))=\sum_{i=0}^e\dim_k(\M^i/\M^{i+1})$. We recall the definition of the {\sl level} of a local, Artinian $k$--algebra (see [Re], Section 5).

\definition{Definition 2.5}
Let $A$ be a local, Artinian $k$--algebra. If $\M^e\ne0$ and $\M^{e+1}=0$ we define the
level of $A$ as $e$ and denote it by $\lev(A)$.

If $e=\lev(A)$ and $n_i:=\dim_k(\M^i/\M^{i+1})$, $0\le i\le e$, we define the Hilbert function of $A$ as the vector $H(A):=(n_0,\dots,n_e)\in{\Bbb N}^{e+1}$. 
\enddefinition

Notice that some other authors prefer to use {\sl maximum socle degree}\/ instead of level.

In any case $n_0=1$. Recall that the Gorenstein condition is equivalent to saying that the socle $\Soc(A):=0\colon \M$ of $A$ is a vector space over $k\cong A/\M$ of dimension $1$. If $e=\lev(A)\ge1$ trivially $\M^e\subseteq\Soc(A)$, hence if $A$ is Gorenstein then equality must hold and $n_e=1$, thus if $\emdim(A)\ge2$ we deduce that $\lev(A)\ge2$ and $\deg(A)\ge\emdim(A)+2$.

In particular, taking into account of Sections 5F.i.a, 5F.i.c and 5F.ii.a of [Ia3] (see also [Ia2]), the list of all possible shapes of Hilbert functions of local, Artinian, Gorenstein $k$--algebra $A$ with $\emdim(A)\ge4$ of degree $7$ is
$$
(1,4,1,1),\quad(1,5,1),
\tag2.6
$$
in degree $d=8$ is
$$
(1,4,1,1,1),\quad(1,5,1,1),\quad(1,6,1),\quad(1,4,2,1),
\tag2.7
$$
in degree $d=9$ is
$$
(1,4,1,1,1,1),\quad(1,5,1,1,1),\quad(1,6,1,1),\quad(1,7,1),\quad(1,4,2,1,1),\quad(1,4,3,1).
\tag2.8
$$
All the sequences in (2.6), (2.7) and (2.8) above actually occur as Hilbert functions of a local, Artinian, Gorenstein $k$--algebra. They can be divided into three different families according to $\dim_k(\M^2/\M^3)$

When $\dim_k(\M^2/\M^3)=1$ the above sequences completely characterize the algebra if $\char(k)\ne2$ (see [Sa]. Another proof can be found in [C--N]: in the proof we need only that $2$ is invertible in $k$), since for a local, Artinian $k$--algebra $A$ of degree $d\ge n+2$, one has $H(A)=(1,n,1,\dots,1)$ if and only if $A\cong A_{n,d}$ where
$$
A_{n,d}:=k[x_1,\dots,x_n]/(x_ix_j, x_h^2-x_1^{d-n},x_1^{d-n+1})_{1\le i<j\le n,\atop 2\le h\le n}.
$$
\remark{Remark 2.9}
In order to describe $\Hilb_d^G(\p N)$, it is interesting to notice that the algebra $A_{n,d}$ is a flat specialization of the simpler algebra $A_{n,d-1}\oplus A_{0,1}$, for each $d\ge n+2\ge4$. Indeed in $k[b,x_1,\dots,x_n]$ we have
$$
\aligned
J:=(&x_ix_j, x_h^2-bx_1^{d-n-1}-x_1^{d-n},x_1^{d-n+1})_{1\le i<j\le n,\atop 2\le h\le n}=\\
&=(x_1+b,x_2,\dots,x_n)\cap (x_ix_j, x_h^2-bx_1^{d-n-1},x_1^{d-n})_{1\le i<j\le n,\atop 2\le h\le n},
\endaligned\tag2.9
$$
for each $d\ge n+2\ge4$. Thus ${\Cal A}_{n,d}:=k[b,x_1,\dots,x_n]/J\to \a1$ is a flat family having special fibre over $b=0$ isomorphic to  $A_{n,d}$ and general fibre isomorphic to $A_{n,d-1}\oplus A_{0,1}$. 
\endremark
\medbreak

In the next two sections we will classify the two remaining cases. More precisely in Section 4we deal with the case $\dim_k(\M^2/\M^3)=2$, i.e. $H(A)=(1,n,2,1,\dots,1)$, $n\ge2$. We obtain the same results proved in [E--V] under the restrictive hypothesis $\char(k)=0$. Finally, in Section 5, we will examine the remaining case, namely $\dim_k(\M^2/\M^3)=3$, i.e. $H(A)=(1,n,3,1)$, $n\ge3$.

\head
3. $k$--Algebras with Hilbert function $(1,n,2,1,\dots,1)$
\endhead
Let $A$ be a local, Artinian, Gorenstein $k$--algebra with $H(A)=(1,n,2,1,\dots,1)$ where $n\ge 2$. The level, $e:=\lev(A)$, is then equal to $d-n-1\ge 3$ and we will assume in this section that $\char(k)>e\ge3$. Consider generators $a_1,\dots,a_n\in\M$.  From $\dim_k(\M^2/\M^3)=2$ and $\M^2=(a_ia_j)_{i,j=1,\dots,n}$, if $a_1a_2\in\M^2\setminus\M^3$, from 
$(a_1+a_2)^2=a_1^2+2a_1a_2+a_2^2$, then at least one among $(a_1+a_2)^2$, $a_1^2$, $a_2^2$ is not in $\M^3$. Thus, up to a linear change of the minimal generators of $\M$, we can always assume that $a_1^2\in\M^2\setminus\M^3$, hence the following three cases are possible for $\M^2$:
$$
(a_1^2,a_1a_2),\qquad (a_1^2,a_2^2),\qquad (a_1^2,a_2a_3).
$$
Let us examine the last case which occurs only if $n\ge3$. As above, since $(a_2+a_3)^2=a_2^2+2a_2a_3+a_3^2$, it follows again, at least one among $(a_2+a_3)^2$, $a_2^2$, $a_3^2$ is not in $(a_1^2)$, i.e. we can finally assume $\M^2=(a_1^2,a_2^2)$ in this case, for a suitable set of minimal generators of $\M$.

It follows the existence of a non--trivial relation of the form $\alpha_1 a_1^2+\alpha_2  a_2^2+\overline{\alpha} a_1a_2\in\M^3$, where $\alpha_1,\alpha_2,\overline{\alpha}\in
k\subseteq A$. The first member of the above relation can be interpreted as the defining polynomial of a single quadric $Q$ in the projective space ${\Bbb P}(V)$ associated to the subspace $V\subseteq\M/\M^2$ generated by the classes of $a_1,a_2$. Such a quadric has rank either $2$ or $1$.

\subhead
3.1. The case $\rk(Q)=2$
\endsubhead
In the second case, again via a suitable linear transformation in $V$ we can assume  $a_1a_2\in\M^3$, thus  $\M^h=(a_1^h,a_2^h)$, for each $h\ge2$. In particular, possibly interchanging $a_1$ and $a_2$, we can assume $a_1^e\ne0$. Thus we obtain both $\M^h=(a_1^h)$, $h\ge3$ and the relations
$$
a_ia_j=\alpha_{i,j}^1a_1^2+\alpha_{i,j}^2a_2^2+\alpha_{i,j}a_1^3,\qquad i,j\ge1,
\tag 3.1.1
$$
where ${\alpha}_{i,j}=\sum_{h=0}^{e-4}\beta_{i,j}^ha_1^h+\beta_{i,j}a_1^{e-3}$, $\alpha_{i,j}^h,\beta_{i,j}^h,\beta_{i,j}\in k$, $\alpha_{i,j}^h=\alpha_{j,i}^h$, $\alpha_{i,j}=\alpha_{j,i}$, $\alpha_{1,1}^1=\alpha_{2,2}^2=1$ and $\alpha_{1,1}^2=\alpha_{1,1}=\alpha_{1,2}^1=\alpha_{1,2}^2=\alpha_{2,2}^1=\alpha_{2,2}=0$.

Via $a_2\mapsto a_2+\alpha_{1,2} a_1^2$ we can assume 
$$
a_1a_2=0,\tag3.1.2
$$
i.e. $\alpha_{1,2}=0$.

Again via $a_j\mapsto a_j+\alpha_{1,j}^1a_1+\alpha_{2,j}^2a_2+\alpha_{1,j}a_1^2$, $j\ge3$, we can assume $\alpha_{1,j}^1=\alpha_{1,j}=0$, $j\ge2$, and $\alpha_{2,j}^2=0$, $j\ne2$. In particular $a_1^2a_j=a_2a_ia_j=0$, $j\ne3$. Explicitely $a_1a_j=\alpha_{1,j}^2a_2^2$, $j\ge2$, and  $a_2a_j=\alpha_{2,j}^2a_2^2+\alpha_{2,j}a_1^3$, $j\ge2$. 

Moreover $a_2^3=\sum_{i=3}^e\mu_ia_1^i$, $\mu_i\in k$, then $\sum_{i=3}^{e-1}\mu_ia_1^{i+1}=a_1a_2^3=0$, thus $\mu_i=0$, $i=3,\dots,e-1$, whence $a_2^3=\mu_ea_1^e$. If $\mu_e=0$ then $a_2^2\in\Soc(A)\setminus\M^e$. Up to multiplying $a_1$ we can thus assume
$$
a_2^3-a_1^e=0.\tag3.1.3
$$
If $n=2$ we have finished so, from now on, we will assume $n\ge3$.

Since $\alpha_{2,j}^{2}=0$, $j\ne 2$, it follows
$\alpha_{i,j}^{2}a_2^{3}=a_2(a_ia_j)=a_2a_ia_j=(a_2a_j)a_i=0$, $(i,j)\ne(2,2)$, whence
$\alpha_{i,j}^2=0$, $(i,j)\ne(2,2)$, thus
$$
a_1a_j=0.\qquad j\ne1. \tag3.1.4
$$
Moreover $\alpha_{i,j}^1a_1^3+\alpha_{i,j}a_1^4=(a_ia_j)a_1=a_1a_ia_j=(a_1a_j)a_i=0$, thus $\alpha_{i,j}^1=0$ and necessarily $\alpha_{i,j}=\beta_{i,j}a_1^{e-3}$. Recall that $\beta_{1,j}=\beta_{2,2}=0$.

Let $y:=y_0+\sum_{i=1}^ny_ia_i+y_{n+1}a_1^2+y_{n+2}a_2^2+\sum_{h=3}^ey_{n+h}a_1^h\in\Soc(A)$, $y_h\in k$. Then the conditions $a_jy=0$, $j=1,\dots,n$, become
$$
\cases
y_0a_1+y_1a_1^2+y_{n+1}a_1^3+\sum_{h=3}^{e-1}y_{n+h}a_1^{h+1}=0,\\
y_0a_2+y_2a_2^2+\left(\sum_{i=3}^{n}y_{i}\beta_{2,i}\right)a_1^e+y_{n+2}a_1^e=0,\\
y_0a_j+y_2\beta_{2,j}a_1^e+\left(\sum_{i=3}^{n}y_{i}\beta_{i,j}\right)a_1^e=0.\qquad j\ge3
\endcases
$$
It is clear that $y_0=y_1=y_2=y_{n+1}=y_{n+3}=\dots=y_{n+e-1}=0$ and $\sum_{i=3}^{n}y_{i}\beta_{i,j}=0$, $j\ge3$. If the symmetric matrix $B:=\left(\beta_{i,j}\right)_{i,j\ge3}$ is singular  then it would be easy to find $y\in\Soc(A)\setminus\M^e$. We conclude that we can make a linear change on $a_3,\dots,a_n$ in such a way that
$$
a_ia_j=\delta_{i,j}a_1^e,\qquad i,j\ge3. \tag3.1.5
$$
Now we finally have $a_2a_j=\gamma_{j}a_1^e$: via $a_2\mapsto a_2+\sum_{j=3}^n\gamma_{j}a_j$ we also obtain
$$
a_2a_j=0,\qquad j\ge3. \tag3.1.6
$$
Combining Equalities (3.1.2), (3.1.3), (3.1.4), (3.1.5) and (3.1.6), we obtain the isomorphism $A\cong A_{n,2,d}^2$ where
$$
A_{n,2,d}^2:=k[x_1,\dots,x_n]/(x_1x_2,x_2^3-x_1^{d-n-1},x_ix_j,x_h^2-x_1^{d-n-1},x_1^{d-n})_{{1\le i< j\le n,\ 3\le j}\atop 3\le h\le n}.
$$
\medbreak

\subhead
3.2. The case $\rk(Q)=1$
\endsubhead
In the first case, via a suitable linear transformation in $V$ we can assume $\alpha_1=\overline{\alpha}=0$, $\alpha_2=1$ in the above relation, i.e. $a_2^2\in\M^3$, thus $\M^h=(a_1^h,a_1^{h-1}a_2)$, $h\ge2$. In particular $\M^e=(a_1^e,a_1^{e-1}a_2)$ and $a_1^{e-t}a_2^t=0$, $t\ge2$. If $a_1^e=0$ then $a_1^{e-1}a_2\ne0$, hence $(a_1+a_2)^e=a_1^e+ea_1^{e-1}a_2\ne0$, hence the linear change $a_1\mapsto a_1+a_2$ allows us to assume both $\M^h=(a_1^h)$, $h\ge3$, and the relations
$$
a_ia_j=\alpha_{i,j}^1a_1^2+\alpha_{i,j}^2a_1a_2+\alpha_{i,j}a_1^3,\qquad i,j\ge1,
\tag 3.2.1
$$
where ${\alpha}_{i,j}=\sum_{h=0}^{e-4}\beta_{i,j}^ha_1^h+\beta_{i,j}a_1^{e-3}$, $\alpha_{i,j}^h,\beta_{i,j}^h,\beta_{i,j}\in k$,  $\alpha_{i,j}^h=\alpha_{j,i}^h$, $\alpha_{i,j}=\alpha_{j,i}$, $\alpha_{1,1}^1=\alpha_{1,2}^2=1$ and $\alpha_{1,1}^2=\alpha_{1,1}=\alpha_{1,2}^1=\alpha_{1,2}=\alpha_{2,2}^1=\alpha_{2,2}^2=0$.

Via the transformation $a_j\mapsto a_j+\alpha_{i,j}^1a_1+\alpha_{i,j}^2a_2+\alpha_{i,j}a_1^2$ we can assume
$$
a_1a_j=0\qquad j\ge3.\tag 3.2.2
$$
Moreover $a_1^2a_2=\sum_{i=3}^e\mu_ia_1^i\in\M^3$, where $\mu_i\in k$.

From now on we will assume $e\ge4$ and we will came back to the case $e=3$ later on. Since $e\ge4$, it follows $a_1^{e-1}a_2=(a_1^2a_2)a_1^{e-3}=\mu_{3}a_1^e$. On the other hand $\mu_3^2a_1^{e}=\mu_3a_1^{e-1}a_2=a_1^{e-4}(a_1^2a_2)a_2=a_1^{e-2}a_2^2=0$,
whence $\mu_3=0$. Via $a_2\mapsto a_2+\sum_{i=4}^e\mu_ia_1^{i-2}$ we obtain
$$
a_1^2a_2=0.\tag 3.2.3
$$
Since $a_2^2\in\M^3$, we have $a_2^2=\sum_{h=3}^{e}\gamma_ha_1^h$. Assume that $\gamma_h=0$, $h<t$. Then $\gamma_ta_1^e=a_1^{e-t}a_2^2=a_1^{e-t-2}(a_1^{2}a_2)a_2=0$. It follows that $\gamma_h=0$, $h\le e-2$, thus $a_2^2=\gamma_{e-1}a_1^{e-1}+\gamma_ea_1^e$. If $\gamma_{e-1}=0$ then $a_1a_2\in\Soc(A)\setminus\M^e$, a contradiction, since $A$ Gorenstein implies $\Soc(A)=\M^e$. Thus $\gamma_{e-1}\ne0$, hence we can find a square root $u$ of $\gamma_{e-1}+\gamma_ea_1$ and via $a_2\mapsto ua_2$ we finally obtain
$$
a_2^2-a_1^{e-1}=0.\tag 3.2.4
$$
If $n=2$ we have finished so, from now on, we will assume $n\ge3$.

Equalities (3.2.2) and (3.2.3) yield $\alpha_{i,j}^1a_1^3+\alpha_{i,j}a_1^4=a_1(a_ia_j)=a_1a_ia_j(a_1a_i)a_j=0$, $(i,j)\ne(1,1),(1,2),(2,2)$, thus $\alpha_{i,j}^1=0$ and $\alpha_{i,j}=\beta_{i,j}^h=0$, $h=0,\dots,e-4$. Moreover Equalities (3.2.2), (3.2.3) and (3.2.4) imply $\alpha_{i,j}^2a_1^e=\alpha_{i,j}^2a_1a_2^2=a_2(a_ia_j)=(a_2a_i)a_j=0$, $(i,j)\ne(1,1),(1,2),(2,2)$, whence $\alpha_{i,j}^2=0$ too. 

Let $y:=y_0+\sum_{i=1}^ny_ia_i+y_{n+1}a_1^2+y_{n+2}a_1a_2+\sum_{h=3}^ey_{n+h}a_1^h\in\Soc(A)$, $y_h\in k$. Then the conditions $a_jy=0$, $j=1,\dots,n$, become
$$
\cases
y_0a_1+y_1a_1^2+y_2a_1a_2+y_{n+1}a_1^3+\sum_{h=3}^{e-1}y_{n+h}a_1^{h+1}=0,\\
y_0a_2+y_1a_1a_2+y_2a_1^{e-1}+\left(\sum_{i=3}^{n}y_{i}\beta_{2,i}\right)a_1^e+y_{n+2}a_1^e=0,\\
y_0a_j+y_2\beta_{2,j}a_1^e+\left(\sum_{i=3}^{n}y_{i}\beta_{i,j}\right)a_1^e=0,\qquad j\ge3.
\endcases
$$
It is clear that $y_0=y_1=y_2=y_{n+1}=y_{n+3}=\dots=y_{n+e-1}=0$ and $\sum_{i=3}^{n}y_{i}\beta_{i,j}=0$, $j\ge3$. If the symmetric matrix $B:=\left(\beta_{i,j}\right)_{i,j\ge3}$ is singular, again $\Soc(A)\ne\M^e$. We conclude that there exists $P\in\GL_{n-3}(k)$ such that ${}^tPBP=I_{n-3}$ is the identity. This matrix corresponds to a linear change of the generators $a_3,\dots,a_n$ which allows us to assume
$$
a_ia_j=\delta_{i,j}a_1^e,\qquad i,j\ge3. \tag3.2.5
$$
At this point we have $a_2a_j=\vartheta_ja_1^3$. Via $a_2\mapsto a_2/2+\sum_{i=3}^{n}y_{i}\vartheta_{i}a_i$ we finally obtain $a_2a_j=0$, $j\ge 3$, and
$a_2^2=a_1^{e-1}+\lambda a_1^e$ for a suitable $\lambda\in A$. Let $v$ be a square root of $1+\lambda a_1$. Then via $a_2\mapsto v a_2$ we finally obtain again Equality (3.2.4) and also
$$
a_2a_j=0,\qquad j\ge3.\tag3.2.6
$$

If $e=3$ then $a_1^2a_2=\mu_3a_1^3$, $a_2^2=\nu_3a_1^3$. If $\mu_3=0$ then $a_1a_2\in\Soc(A)\setminus\M^e$, thus up to multiplying $a_2$ by a suitable square root of $\mu_e$, we can assume $\mu_e=1$. Via $a_2\mapsto a_2+\beta_{2,2}a_1^2/2$ we finally obtain
$$
a_1^2a_2-a_1^3=a_2^2=0.\tag 3.2.7
$$
If $n=2$ we have finished. Otherwise, if $n\ge3$, we can repeat word by word the discussion above and we finally obtain Equalities (3.2.5) and (3.2.6) with $\lambda=0$.

In particular, combining Equalities (3.2.2), (3.2.3), (3.2.4), (3.2.5), (3.2.6) and (3.2.7), we obtain the isomorphism $A\cong A_{n,2,d}^1:=k[x_1,\dots,x_n]/I$ where
$$
I=\cases
(x_1^2x_2-x_1^3,x_2^2,x_ix_j,x_h^2-x_1^{3},x_1^{4})_{{1\le i<j\le n,\ 3\le j}\atop 3\le h\le n}
&\text{if $d=n+4$,}\\(x_1^2x_2,x_2^2-x_1^{d-n-2},x_ix_j,x_h^2-x_1^{d-n-1},x_1^{d-n})_{{1\le i<j\le n,\ 3\le j}\atop 3\le h\le n}
&\text{if $d\ge n+5$.}
\endcases
$$

It is natural to investigate if we have determined non--isomorphic algebras (see Proposition 4.5 of [C--N]).

\proclaim{Proposition 3.3}
The algebras $A_{n,2,d}^{1}$ and $A_{n,2,d}^{2}$ are non--isomorphic.
\endproclaim
\demo{Proof}
Set $I_h$ be the ideal generated by $\{\ u\in A_{n,2,d}^{h}\ \vert\ u^2\in\M^e\ \}$.
Then $\sum_{i=1}^n\lambda_i\overline{x_i}\in I_h$ if and only if 
$$
\sum_{i,j=1}^n\lambda_i\lambda_j\overline{x}_i \overline{x}_j=\lambda_1^2\overline{x}_1^2+\lambda_1\lambda_2\overline{x}_1 \overline{x}_2+\lambda_2^2\overline{x}_2^2+\sum_{i=3}^n\lambda_i^2\overline{x}_i^2\in \M^e.
$$
This is equivalent to $\lambda_i=0$, $i\le h$. In particular $I_h=(\overline{x}_{h+1},\dots,\overline{x}_n)$, hence $\dim_k(I_h\otimes k)=n-h$.

If $A_{n,2,d}^{1}\cong A_{n,2,d}^{2}$ then the ideals $I_1$ and $I_2$ would correspond each other in this isomorphism, hence such dimension should coincide.
\qed
\enddemo

\remark{Remark 3.4}
As showed for $A_{n,d}$, also the algebras $A_{n,2,d}^{h}$ are flat specialization of easier algebras. Indeed take for $h=2$
$$
\align
J&:=(x_1x_2,x_2^3-bx_1^{d-n-2}-x_1^{d-n-1},x_ix_j,x_h^2-bx_1^{d-n-2}-x_1^{d-n-1},x_1^{d-n})_{{1\le i\le n,\atop 3\le j\le n,i\ne j}\atop 3\le h\le n}=\\
&=(x_1+b,x_2,\dots,x_n)\cap (x_1x_2,x_2^3-bx_1^{d-n-2},x_ix_j,x_h^2-bx_1^{d-n-2},x_1^{d-n-1})_{{1\le i\le n,\atop 3\le j\le n,i\ne j}\atop 3\le h\le n}.
\endalign
$$
Thus ${\Cal A}:=k[b,x_1,\dots,x_n]/J\to \a1$ is a flat family having special fibre over $b=0$ isomorphic to  $A_{n,2,d}^{2}$ and general fibre isomorphic to $A_{n,2,d-1}^{2}\oplus A_{0,1}$ if $e\ge 4$ and $A_{n,n+3}\oplus A_{0,1}$ if $d-n-1=e=3$.

Finally consider the case of the $k$--algebras $A_{n,2,d}^{1}$. Then let us consider the ideal $J$ defined as
$$
(bx_1x_2+x_1^2,x_1^2x_2+bx_2^3-x_1^3,x_ix_j,x_h^2-x_1^{3},x_1^{4})_{{1\le i\le n,\atop 3\le j\le n,i\ne j}\atop 3\le h\le n}
$$
if  $d=n+4$ (i.e. $e=3$), and
$$
(bx_1x_2+x_2^2-x_1^{d-n-2},bx_2^3-bx_1^{d-n-1}+x_1^2x_2,x_ix_j,x_h^2-x_1^{d-n-1},x_1^{d-n})_{{1\le i\le n,\atop 3\le j\le n,i\ne j}\atop 3\le h\le n}
$$
if $d\ge n+5$ (i.e. $e\ge4$).
In this case ${\Cal A}:=k[b,x_1,\dots,x_n]/J\to \a1$ is a flat family of local Artinian, Gorenstein $k$--algebras with constant Hilbert function $H({\Cal A}_b)=(1,n,2,1,\dots,1)$. For $b=0$ we have trivially ${\Cal A}_0\cong A_{n,2,d}^{1}$. For general $b\ne0$ the algebra ${\Cal A}_b$ is again local and Gorenstein, thus it must be either $A_{n,2,d}^{1}$ or $A_{n,2,d}^{2}$.

In any case we have the relations $bx_1x_2+x_1^2\in\M^3$ (if $d=n+4$) or $bx_1x_2+x_2^2\in\M^3$ (if $d\ge n+5$) in ${\Cal A}_b$, thus, computing the invariant $\dim_k(I_h\otimes k)$ defined in the proof of Proposition 3.3 in these cases, we finally obtain ${\Cal A}_b\cong A_{n,2,d}^{2}$ for general $b\ne0$.
\endremark
\medbreak

\head
4. $k$--Algebras with Hilbert functions $(1,n,3,1)$
\endhead
Let $\M=(a_1,a_2,a_3,\dots,a_n)$. As in the previous section, one can always assume that $\M^2=(a_1^2,a_2^2,a_3^2)$. Thus we have three linearly independent relations of the form
$$
\alpha_1 a_1^2+\alpha_2  a_2^2+\alpha_3a_3^2+2\overline{\alpha}_1 a_2a_3+2\overline{\alpha}_2 a_1a_3+2\overline{\alpha}_3 a_1a_2\in\M^3,\tag 4.1
$$
where $\alpha_i,\overline{\alpha}_j\in
k\subseteq A$, $i,j=1,2,3$, and hence a net $\Cal N$ of conics in the projective space ${\Bbb P}(V)$, associated to the subspace $V\subseteq \M/\M^2$ generated by the classes of $a_1,a_2,a_3$. Let $\Delta$ be the discriminant curve of $\Cal N$ in ${\Bbb P}(V)$. Then $\Delta$ is a plane cubic and the classification of $\Cal N$ depends on the structure of $\Delta$ as explained in [Wa].

We stress that the classification described in [Wa] in the case complex case, holds as soon as the base field $k$ is an algebraically closed field with $\char(k)\ne2,3$.

\subhead
4.2  The case of integral $\Delta$
\endsubhead
Let $\Delta$ be irreducible. 
Then, taking into account the results proved in [Wa], we obtain that Relations (4.1) above become
$a_1a_2+a_3^2,a_1a_3,a_2^2-6p a_3^2+q a_1^2\in\M^3$, where $p,q\in k$. In particular
$$
\gathered
a_1^2a_2=a_1^2a_3=a_1a_2a_3=a_1a_3^2=a_2^2a_3=a_3^3=0,\\
a_2a_3^2=-a_1a_2^2=q a_1^3,\quad
a_2^3=-6pq a_1^3,
\endgathered\tag4.2.1
$$
thus
$\M^2=(a_1^2,a_3^2,a_2a_3)$ and $\M^3=(a_1^3)$. Relation (4.2) thus become
$$
a_1a_2=-a_3^2+\beta_{1,2} a_1^3,\qquad a_1a_3=\beta_{1,3} a_1^3,\qquad a_2^2=\alpha_{2,2}^1 a_1^2+\alpha_{2,2}^3a_3^2+\beta_{2,2} a_1^3,
$$
where $\alpha_{i,j}^h,{\beta}_{i,j}\in k$, $\alpha_{2,2}^1=-q$, $\alpha_{2,2}^3=6p$..

In general, we have relations of the form
$$
a_ia_j=\alpha_{i,j}^1a_1^2+\alpha_{i,j}^2a_2a_3+\alpha_{i,j}^3a_3^2+\beta_{i,j}a_1^3,\qquad i\ge1,j\ge4,\tag4.2.2
$$
where the ${\alpha}_{i,j}^h,\beta_{i,j}\in k$ are as above and $\alpha_{i,j}^h=\alpha_{j,i}^h$, $\beta_{i,j}=\beta_{j,i}$.

Via $(a_2,a_3)\mapsto(a_2+\beta_{1,2} a_1^2,a_3+\beta_{1,3} a_1^2)$, we can assume $\beta_{1,2}=\beta_{1,3}=0$.

If $\alpha_{2,2}^1=0$ then $a_2a_3\in\Soc(A)\setminus\M^{3}$, a contradiction. Let $u$ be a fourth root of $\beta_{2,2}a_1-\alpha_{2,2}^1$. Via $(a_2,a_3)\mapsto(u^2 a_2,u a_3)$ then we can assume $\alpha_{2,2}^1=-1$ and $\beta_{2,2}=0$, whence
$$
a_1a_2=-a_3^2,\qquad
a_1a_3=0,\qquad
a_2^2=-a_1^2+\alpha a_3^2,\tag 4.2.3
$$
where $\alpha:=\alpha_{2,2}^3$.

If $n=3$ we have finished, thus we will assume $n\ge4$ from now on. Via $a_j\mapsto a_j+\alpha_{1,j}^1a_1+\beta_{1,j}a_1^2+\alpha_{2,j}^2a_3+\beta_{2,j}a_3^2+\alpha_{3,j}^2a_2+\beta_{3,j}a_2a_3$, we can assume $\alpha_{1,j}^1=\beta_{1,j}=\alpha_{2,j}^2=\beta_{2,j}=\alpha_{3,j}^2=\beta_{3,j}=0$, $j\ge4$.

Since $a_1a_3=0$, it follows $\alpha_{1,j}^2a_2a_3^2=(a_1a_j)a_3=a_1a_3a_j=(a_3a_j)a_1=\alpha_{3,j}^1a_1^3$, thus $\alpha_{1,j}^2=\alpha_{3,j}^1=0$, $j\ge4$.
Since $a_1a_2a_j=-a_3^2a_j=(a_3a_j)a_3=0$, we have $\alpha_{1,j}^3a_2a_3^2=(a_1a_j)a_2=a_1a_2a_j=(a_2a_j)a_1=\alpha_{2,j}^1a_1^3$, thus $\alpha_{1,j}^3=\alpha_{2,j}^1=0$, $j\ge4$.
Moreover $0=(a_2a_j)a_3=a_2a_3a_j=(a_3a_j)a_2=\alpha_{3,j}^3a_2a_3^2$, thus $\alpha_{3,j}^3=0$, $j\ge4$.
Finally $0=a_2^2a_j=(a_2a_j)a_2=\alpha_{2,j}^3a_2a_3^2$, thus $\alpha_{3,j}^3=0$, $j\ge4$. We conclude that $a_1a_j=a_2a_j=a_3a_j=0$, $j\ge4$.

It follows that $0=(a_1a_i)a_j=(a_ia_j)a_1=\alpha_{i,j}^1a_1^3$, $0=(a_2a_i)a_j=(a_ia_j)a_2=\alpha_{i,j}^3a_2a_3^2$, $0=(a_3a_i)a_j=(a_ia_j)a_3=\alpha_{i,j}^2a_2a_3^2$, thus $a_ia_j=\beta_{i,j}a_1^3$, $i,j\ge4$.

Let $y:=y_0+\sum_{i=1}^ny_ia_i+y_{n+1}a_1^2+y_{n+2}a_2a_3+y_{n+3}a_3^2+y_{n+4}a_1^{3}\in\Soc(A)$, $y_h\in k$. Then the conditions $a_jy=0$, $j=1,\dots,n$, become
$$
\cases
y_0a_1+y_1a_1^2-y_2a_3^2+y_{n+1}a_1^3=0,\\
y_0a_2-y_1a_3^2+y_2\alpha a_3^2+y_2a_1^{e-1}+y_3a_2a_3+y_{n+3}a_1^3=0,\\
y_0a_3+y_2a_2a_3+y_3a_3^2+y_{n+2}a_1^3=0,\\
y_0a_j+\left(\sum_{i=4}^{n}y_{i}\beta_{i,j}\right)a_1^3=0,\qquad j\ge4.
\endcases
$$
It is clear that $y_0=y_1=y_2=y_3=y_{n+1}=y_{n+3}=0$ and $\sum_{i=4}^{n}y_{i}\beta_{i,j}=0$, $j\ge4$. If the symmetric matrix $B:=\left(\beta_{i,j}\right)_{i,j\ge4}$ would be singular again $\Soc(A)\ne\M^e$. We conclude that we can make a linear change on $a_4,\dots,a_n$ in such a way that
$$
a_ia_j=\delta_{i,j}a_1^3,\qquad i,j\ge4. \tag4.2.4
$$
Combining Equalities (4.2.2), (4.2.3) and (4.2.4) we obtain
$$
A\cong A_{n,3,n+5}^{1,\alpha^2}:=k[x_1,\dots,x_n]/(x_1x_2+x_3^2,x_1x_3,x_2^2-\alpha x_3^2+x_1^{2},x_ix_j,x_h^2-x_1^3)_{{1\le i<j\le n,\ 4\le j}\atop 4\le h\le n}.
$$

Let $i^2=-1$: we finally notice that via $(a_2,a_3)\mapsto(-a_2, ia_3)$ we can identify the two cases $\pm\alpha$.

\subhead
4.3. The case of non--integral $\Delta$
\endsubhead
In this case we have an easy classification described in [Wa] for the possible Relations (4.2).

\subhead
4.3.1. The cases $D$ and $E$
\endsubhead
In case $D$ Relations (4.1) become
$a_1^2,a_2^2, a_3^2+2pa_1a_2\in\M^3$ where $p=1$ in case $D$ and $p=0$ in case $E$. In particular
$$
a_1^2a_2=a_1a_2^2=a_1^2a_3=a_1a_3^2=a_2^2a_3=a_2a_3^2=0,
\tag4.3.1.1
$$
thus
$\M^2=(a_1a_2,a_1a_3,a_2a_3)$ and $\M^3=(a_1a_2a_3)$ and Relation (4.1) become
$$
a_1^2=\beta_{1,1} a_1a_2a_3,\qquad a_2^2=\beta_{2,2} a_1a_2a_3,\qquad a_3^2=-2pa_1a_2+\beta_{2,2} a_1a_2a_3,
$$
where ${\beta}_{i,j}\in k$.

Via $(a_1,a_2,a_3)\mapsto (a_1+\beta_{1,1} a_2a_3/2,a_2+\beta_{2,2} a_1a_3/2,a_3+\beta_{3,3} a_1a_2/2)$ we can assume $\beta_{1,1}=\beta_{2,2} =\beta_{3,3}=0$.
In general, we have relations of the form
$$
a_ia_j=\alpha_{i,j}^1a_1a_2+\alpha_{i,j}^2a_1a_3+\alpha_{i,j}^3a_2a_3+\beta_{i,j}a_1a_2a_3,\qquad i\ge1,j\ge4,\tag4.3.1.2
$$
where ${\beta}_{i,j},\alpha_{i,j}^h,\in k$,  $\alpha_{i,j}^h=\alpha_{j,i}^h$, $\beta_{i,j}=\beta_{j,i}$. Via 
$a_j\mapsto a_j+\alpha_{1,j}^1a_2+\alpha_{1,j}^2a_3+\alpha_{2,j}^1a_1+\beta_{1,j}a_2a_3+\beta_{2,j}a_1a_3+\beta_{3,j}a_1a_2$ we can assume also that $\alpha_{1,j}^1=\alpha_{1,j}^2=\alpha_{2,j}^1=\beta_{1,j}=\beta_{2,j}=\beta_{3,j}=0$, $j\ge4$.

Since $a_1^2=0$, we have $0=(a_1^2)a_j=(a_1a_j)a_1=\alpha_{1,j}^3a_1a_2a_3$, hence $\alpha_{1,j}^3=0$, $j\ge4$. Similarly, since $a_2^2=0$ we also obtain $\alpha_{2,j}^2=0$, $j\ge4$. Since $a_1a_j=0$, we have $0=(a_1a_j)a_h=(a_ha_j)a_1=\alpha_{h,j}^3a_1a_2a_3$, $h=2,3$, hence $\alpha_{2,j}^3=\alpha_{3,j}^3=0$, $j\ge4$. Similarly, looking at $a_2a_3a_j$, we also infer $\alpha_{3,j}^2=0$. Finally $0=-2(a_2a_j)a_1=(-2pa_1a_2)a_j=a_3^2a_j=\alpha_{3,j}^1a_1a_2a_3$, thus $\alpha_{3,j}^1=0$, $j\ge4$. We conclude that $a_1a_j=a_2a_j=a_3a_j=0$, $j\ge4$.

As in the previous case, it follows that $0=(a_1a_i)a_j=(a_ia_j)a_1=\alpha_{i,j}^3a_1a_2a_3$, $0=(a_2a_i)a_j=(a_ia_j)a_2=\alpha_{i,j}^2a_1a_2a_3$, $0=(a_3a_i)a_j=(a_ia_j)a_3=\alpha_{i,j}^1a_1a_2a_3$, thus $a_ia_j=\beta_{i,j}a_1a_2a_3$, $i,j\ge4$.

Let $y:=y_0+\sum_{i=1}^ny_ia_i+y_{n+1}a_1a_2+y_{n+2}a_1a_3+y_{n+3}a_2a_3+y_{n+4}a_1a_2a_3\in\Soc(A)$, $y_h\in k$. Then the conditions $a_jy=0$, $j=1,\dots,n$, become
$$
\cases
y_0a_1+y_2a_1a_2+y_3a_1a_3+y_{n+3}a_1a_2a_3=0,\\
y_0a_2+y_1a_1a_2+y_3a_2a_3+y_{n+2}a_1a_2a_3=0,\\
y_0a_3+y_1a_1a_3+y_2a_2a_3-y_3pa_1a_2+y_{n+1}a_1a_2a_3=0,\\
y_0a_j+\left(\sum_{i=4}^{n}y_{i}\beta_{i,j}\right)a_1a_2a_3=0,\qquad j\ge4.
\endcases
$$
It is clear that $y_0=y_1=y_2=y_3=y_{n+3}=0$ and $\sum_{i=4}^{n}y_{i}\beta_{i,j}=0$, $j\ge4$. Again, if the symmetric matrix $B:=\left(\beta_{i,j}\right)_{i,j\ge4}$ would be singular, then $\Soc(A)\ne\M^e$. We conclude that we can make a linear change on $a_4,\dots,a_n$ in such a way that
$$
a_ia_j=\delta_{i,j}a_1a_2a_3,\qquad i,j\ge4. \tag4.3.1.3
$$
Combining Equalities (4.3.1.1), (4.3.1.2) and (4.3.1.3) we get
$$
A\cong A_{n,3,d}^{3-p,0}:=k[x_1,\dots,x_n]/(x_1^2,x_2^2,x_3^2+2px_1x_2,x_ix_j,x_h^2-x_1x_2x_3)_{{1\le i<j\le n,\ 4\le j}\atop 4\le h\le n}.
$$
\medbreak

\subhead
4.3.2. The case $E^*$
\endsubhead
In this case our Relations (4.1) become
$a_1a_3,a_2a_3,a_3^2+2 a_1a_2\in\M^3$. In particular $\M^2=(a_1^2,a_2^2,a_3^2)$ and we have
$$
a_1^2a_3=a_1a_2a_3=a_1a_3^2=a_2^2a_3=a_2a_3^2=a_1^2a_2=a_1a_2^2=a_3^3=0,
$$
thus we can always assume $\M^3=(a_1^3)$, whence
$$
a_1a_3=\beta_{1,3} a_1^3,\qquad a_2a_3=\beta_{2,3} a_1^3,\qquad a_1a_2=\beta_{1,2} a^3_1.
$$
where ${\beta}_{i,j}\in k$.

In general, we have relations of the form
$$
a_ia_j=\alpha_{i,j}^1a_1^2+\alpha_{i,j}^2a_2^2+\alpha_{i,j}^3a_3^2+\beta_{i,j}a_1^3,\qquad i\ge1,j\ge4,
$$
where ${\beta}_{i,j},\alpha_{i,j}^h\in k$,  $\beta_{i,j}=\beta_{j,i}$, $\alpha_{i,j}^h=\alpha_{j,i}^h$.

Let $a_h^3=\lambda_ha_1^3$ for some $\lambda_h\in k$, $h=2,3$. If $\mu_2=0$, then $a_2^2\in\Soc(A)\setminus\M^e$: similarly if $\mu_3=0$. It follows that we can always assume $\mu_2=\mu_3=1$. Thus via $a_3\mapsto a_3+\beta_{2,3}a_2^2$ we also have $\beta_{2,3}=0$. Hence we have
$$
a_1a_2=0,\qquad a_1a_3=0,\qquad a_2a_3=0,\qquad a_2^3=a_1^3,\qquad a_3^3=a^3_1.\tag 4.3.2.1
$$

Since $0=(a_1a_h)a_j=(a_ha_j)a_1=\alpha_{h,j}^1a_1^3$, $h=2,3$, we obtain $\alpha_{2,j}^1=\alpha_{3,j}^1=0$. Thus via 
$a_j\mapsto a_j+\alpha_{1,j}^1a_1+\alpha_{2,j}^2a_2+\alpha_{3,j}^3a_3+\beta_{1,j}a_1^2+\beta_{2,j}a_2^2+\beta_{3,j}a_3^2$ we can assume $\alpha_{1,j}^1=\alpha_{2,j}^2=\alpha_{3,j}^3=\beta_{1,j}=\beta_{2,j}=\beta_{3,j}=0$.

Since $0=(a_1a_h)a_j=(a_1a_j)a_h=\alpha_{1,j}^ha_1^3$, we have $\alpha_{1,j}^h=0$, $h=2,3$. Similarly, since $a_2a_3a_j=0$, one also obtain $\alpha_{2,j}^3=\alpha_{3,j}^2=0$, thus $a_ia_j=0$, $i=1,2,3$, $j\ge4$. It follows that $0=(a_1a_i)a_j=(a_ia_j)a_1=\alpha_{i,j}^1a_1^3$, $0=(a_ha_i)a_j=(a_ia_j)a_h=\alpha_{i,j}^ha_1^3$, $h=2,3$, thus $a_ia_j=\beta_{i,j}a_1^3$, $i,j\ge4$.

Let $y:=y_0+\sum_{i=1}^ny_ia_i+y_{n+1}a_1^2+y_{n+2}a_2^2+y_{n+3}a_3^2+y_{n+4}a_1^{3}\in\Soc(A)$, $y_h\in k$. Then the conditions $a_jy=0$, $j=1,\dots,n$, become
$$
\cases
y_0a_1+y_1a_1^2+y_{n+1}a_1^3=0,\\
y_0a_2+y_2 a_2^2+y_{n+2}a_1^{3}=0,\\
y_0a_3+y_3a_3^2+y_{n+3}a_1^3=0,\\
y_0a_j+\left(\sum_{i=4}^{n}y_{i}\beta_{i,j}\right)a_1^3=0,\qquad j\ge4.
\endcases
$$
and we deduce as in the previous cases that we can make a linear change on $a_4,\dots,a_n$ in such a way that
$$
a_ia_j=\delta_{i,j}a_1^3,\qquad i,j\ge4. \tag4.3.2.2
$$
Combining Equalities (4.3.2.1) and (4.3.2.2) we obtain
$$
A\cong A_{n,3,n+5}^{4,0}:=k[x_1,\dots,x_n]/(x_2^3-x_1^3,x_3^3-x_1^3,x_ix_j,x_h^2-x_1^3)_{1\le i<j\le n\atop 4\le h\le n}.
$$
When $n=3$, the well--known structure theorem for Gorenstein local rings, proved in [B--E], guarantees that $(x_ix_j,x_h^2-x_1^3)_{1\le i<j\le 3\atop 2\le h\le 3}$ is minimally generated by the $4\times4$ pfaffians of a suitable $5\times5$ skew--symmetric matrix $M$. E.g. one may take
$$
M:=\pmatrix
0&0&0&x_1&x_2\\
0&0&x_3&-x_1&0\\
0&-x_3&0&x_2^2&x_1^2\\
-x_1&x_1&-x_2^2&0&-x_3^2\\
-x_2&0&-x_1^2&x_3^2&0
\endpmatrix.
$$
\medbreak

\subhead
4.3.3. The case $G^*$
\endsubhead
In this case our Relations (4.1) become
$a_1^2,a_1a_2,a_2a_3\in\M^3$. In particular $\M^2=(a_3^2,a_2^2,a_1a_3)$ and we have
$$
a_1^3=a_1^2a_2=a_1^2a_3=a_1a_2a_3=a_1a_3^2=a_2^2a_3=a_2a_3^2=0.
$$
We also have relations
$$
a_ia_j=\alpha_{i,j}^1a_3^2+\alpha_{i,j}^2a_2^2+\alpha_{i,j}^3a_1a_3+\gamma_{i,j},\qquad i\ge1,j\ge4,
$$
where ${\gamma}_{i,j}\in\M^3$, $\alpha_{i,j}^h\in k$,  $\alpha_{i,j}=\alpha_{j,i}$,  $\alpha_{i,j}^h=\alpha_{j,i}^h$.

If $a_1a_3^2=a_3^3=0$ then $a_3^2\in\Soc(A)\setminus\M^3$. If $a_3^3=0$, then via $a_3\mapsto a_3+\lambda a_1$, we obtain $a_3^3$. Thus we can always assume $\M^3=(a_3^3)$, whence $\gamma_{i,j}=\beta_{i,j}a_3^3$, where $\beta_{i,j}\in k$. In particular
$$
a_1^2=\beta_{1,1} a_3^3,\qquad a_1a_2=\beta_{1,2} a_3^3,\qquad a_2a_3=\beta_{2,3} a^3_3.\tag4.3.3.1
$$

Let $a_2^3=\mu_1a_3^3$, $a_1a_3^2=\mu_2a_3^3$. If $\mu_1=0$, then $a_2^2\in\Soc(A)\setminus\M^3$. If $\mu_2=0$, then $a_1a_3\in\Soc(A)\setminus\M^3$. Thus we can assume $\mu_1=\mu_2=1$. It follows that, via $(a_1,a_2)\mapsto(a_1+\beta_{1,1}a_3^2/2+\beta_{1,2}a_2^2/2,a_2+\beta_{2,3}a_3^2)$, we obtain $\beta_{1,1}=\beta_{1,3}=\beta_{2,3}=0$.

Since $0=(a_1a_h)a_j=(a_ha_j)a_1=\alpha_{h,j}^1a_1^3$, $h=2,3$, we obtain $\alpha_{2,j}^1=\alpha_{3,j}^1=0$. Thus via 
$a_j\mapsto a_j+\alpha_{1,j}^1a_1+\alpha_{2,j}^2a_2+\alpha_{3,j}^3a_3+\beta_{1,j}a_1^2+\beta_{2,j}a_2^2+\beta_{3,j}a_3^2$ we can assume $\alpha_{1,j}^1=\alpha_{2,j}^2=\alpha_{3,j}^3=\beta_{1,j}=\beta_{2,j}=\beta_{3,j}=0$.

Since $0=(a_1a_h)a_j=(a_1a_j)a_h=\alpha_{1,j}^ha_3^3$, it follows $\alpha_{1,j}^h=0$, $h=1,2$. Similarly $0=(a_1a_2)a_j=(a_2a_j)a_1=\alpha_{2,j}^1a_3^3$, then $\alpha_{2,j}^1=0$. Since $a_2a_3a_j=0$, one also obtain $\alpha_{2,j}^3=\alpha_{3,j}^2=0$. Finally $\alpha_{1,j}^3a_3^3=(a_1a_4)a_3=(a_3a_4)a_1=\alpha_{3,4}^1a_3^3$, whence $\alpha_{1,j}^3=\alpha_{3,4}^1$.

Via $a_j\mapsto a_j+\alpha_{1,j}^3a_3+\beta_{1,j}a_3^2+\alpha_{2,j}^2a_2+\beta_{2,j}a_2^2+\alpha_{3,j}^3a_1+(\beta_{3,j}-\beta_{1,j})a_1a_3)$, we finally obtain $a_ia_j=0$, $i=1,2,3$, $j\ge4$. It follows that $0=(a_ha_i)a_j=(a_ia_j)a_h=\alpha_{i,j}^ha_3^3$, thus $a_ia_j=\beta_{i,j}a_3^3$, $i,j\ge4$.

Let $y:=y_0+\sum_{i=1}^ny_ia_i+y_{n+1}a_3^2+y_{n+2}a_2^2+y_{n+3}a_1a_3+y_{n+4}a_3^{3}\in\Soc(A)$, $y_h\in k$. Then the conditions $a_jy=0$, $j=1,\dots,n$, become
$$
\cases
y_0a_1+y_3a_1a_3+y_{n+1}a_3^3=0,\\
y_0a_2+y_2a_2^2+y_{n+2}a_3^{3}=0,\\
y_0a_3+y_1a_1a_3+y_3a_3^2+y_{n+1}a_3^3+y_{n+3}a_3^3=0,\\
y_0a_j+\left(\sum_{i=4}^{n}y_{i}\beta_{i,j}\right)a_3^3=0,\qquad j\ge4.
\endcases
$$
and we deduce as in the previous cases that we can make a linear change on $a_4,\dots,a_n$ in such a way that
$$
a_ia_j=\delta_{i,j}a_3^3,\qquad i,j\ge4. \tag4.3.3.2
$$
Combining Equalities (4.3.3.1) and (4.3.3.2) we obtain
$$
A\cong A_{n,3,n+5}^{5,0}:=k[x_1,\dots,x_n]/(x_1^2,x_1x_2,x_2x_3,x_2^3-x_3^3,x_1x_3^2-x_3^3,x_ix_j,x_h^2-x_3^3)_{{1\le i<j\le n,\ 4\le j}\atop 4\le h\le n}.
$$
When $n=3$, the ideal defining $A_{n,3,8}^{5,0}$ is generated by the $4\times4$ pfaffians of 
$$
M:=\pmatrix
0&0&x_2&-x_3&x_1\\
0&0&-x_2&x_1&0\\
-x_2&x_2&0&x_3^2&-x_3^2\\
x_3&-x_1&-x_3^2&0&x_2^2\\
-x_1&0&x_3^2&-x_2^2&0
\endpmatrix.
$$

\subhead
4.3.4. The case $H$
\endsubhead
In this case our Relations (4.1) become
$a_1^2,a_1a_2,2a_1a_3+a_2^2\in\M^3$. In particular $\M^2=(a_2^2,a_3^2,a_2a_3)$ and we have
$$
a_1^3=a_1^2a_2=a_1^2a_3=a_1a_2^2=a_2^3=0.
$$
Moreover $a_2^2a_3=2a_1a_3^2$. 

If $a_1a_3^2=0$, then $a_1a_3\in\Soc(A)\setminus\M^3$. Thus we can assume  $\M^3=(a_1a_3^2)$, whence $$
a_1^2=\beta_{1,1} a_1a_3^2,\qquad a_1a_2=\beta_{1,2} a_1a_3^2,\qquad a_1a_3=-a_2^2/2+\beta_{1,3} a_1a_3^2,\tag4.3.4.1
$$
where $\beta_{i,j}\in k$. Via $(a_1,a_2,a_3)\mapsto(a_1+\beta_{1,1}a_3^2/2,a_2+\beta_{1,2}a_3^2,a_3+\beta_{1,3}a_3^2)$, we obtain $\beta_{1,1}=\beta_{1,2}=\beta_{1,3}=0$ in Relations (4.3.4.1).

Let $a_2a_3^2=\mu_1a_1a_3^2$, $a_3^3=\mu_2a_1a_3^2$. The transformation $(a_2,a_3)\mapsto (a_2+\mu_1a_1,a_3+\mu_2a_1/2)$ does not affect Relations (4.3.4.1), thus we can assume $\mu_1=\mu_2=0$.

We also have relations
$$
a_ia_j=\alpha_{i,j}^1a_2^2+\alpha_{i,j}^2a_3^2+\alpha_{i,j}^3a_2a_3+\beta_{i,j}a_1a_3^2,\qquad i\ge1,j\ge4,
$$
where ${\beta}_{i,j},\alpha_{i,j}^h\in k$, $\beta_{i,j}=\beta_{j,i}\in k$,  $\alpha_{i,j}^h=\alpha_{j,i}^h$. Via $a_j\mapsto a_j-2\alpha_{1,j}^1a_1+\beta_{1,j}a_3^2+\alpha_{2,j}^3a_3-\beta_{2,j}a_2a_3/2-2\alpha_{3,j}^1a_1+\alpha_{3,j}^3a_2+\beta_{3,j}a_1a_3$, we can assume $\alpha_{1,j}^1=\beta_{1,j}=\alpha_{2,j}^3=\beta_{2,j}=\alpha_{3,j}^1=\alpha_{3,j}^3=\beta_{3,j}=0$.

Since $0=(a_1a_h)a_j=(a_ha_j)a_1=\alpha_{h,j}^2a_1a_3^2$, $h=1,2,3$, we obtain $\alpha_{1,j}^2=\alpha_{2,j}^2=\alpha_{3,j}^2=0$. Since $0=(a_1a_2)a_j=(a_1a_j)a_2=-2\alpha_{1,j}^3a_1a_3^2$, it follows $\alpha_{1,j}^3=0$. Since $0=(a_2a_3)a_j=(a_3a_j)a_2=-2\alpha_{2,j}^1a_1a_3^2$, we obtain $\alpha_{2,j}^1=0$. Thus $a_ia_j=0$, $i=1,2,3$, $j\ge4$.

It follows that $0=(a_1a_i)a_j=(a_ia_j)a_1=\alpha_{i,j}^2a_1a_3^2$, $0=(a_2a_i)a_j=(a_ia_j)a_2=-2\alpha_{i,j}^3a_1a_3^2$, , $0=(a_3a_i)a_j=(a_ia_j)a_3=-2\alpha_{i,j}^1a_1a_3^2$, thus $a_ia_j=\beta_{i,j}a_1a_3^2$, $i,j\ge4$.

Again let $y:=y_0+\sum_{i=1}^ny_ia_i+y_{n+1}a_2^2+y_{n+2}a_3^2+y_{n+3}a_2a_3+y_{n+4}a_1a_3^{2}\in\Soc(A)$, $y_h\in k$. Then the conditions $a_jy=0$, $j=1,\dots,n$, become
$$
\cases
y_0a_1-y_3a_2^2/2+y_{n+2}a_1a_3^2=0,\\
y_0a_2+y_2a_2^2+y_3a_2a_3-2y_{n+3}a_1a_3^{2}=0,\\
y_0a_3-2y_1a_2^2/2+y_2a_2a_3+y_3a_3^2-2y_{n+1}a_1a_3^2=0,\\
y_0a_j+\left(\sum_{i=4}^{n}y_{i}\beta_{i,j}\right)a_1a_3^2=0,\qquad j\ge4
\endcases
$$
and we deduce as in the previous cases that we can make a linear change on $a_4,\dots,a_n$ in such a way that
$$
a_ia_j=\delta_{i,j}a_1a_3^2,\qquad i,j\ge4. \tag4.3.4.2
$$
Combining Equalities (4.3.4.1) and (4.3.4.2) we obtain
$$
A\cong A_{n,3,n+5}^{6,0}:=k[x_1,\dots,x_n]/(x_1^2,x_1x_2,2x_1x_3+x_2^2,x_3^3,x_2x_3^2,x_ix_j,x_h^2-x_1x_3^2)_{{1\le i<j\le n,\ 4\le j}\atop 4\le h\le n}.
$$
When $n=3$, the ideal defining $A_{n,3,8}^{6,0}$ is generated by the $4\times4$ pfaffians of 
$$
M:=\pmatrix
0&0&-2x_3&x_1&-x_2\\
0&0&-x_2&0&x_1\\
2x_3&x_2&0&0&0\\
-x_1&0&0&0&x_3^2\\
x_2&-x_1&0&-x_3^2&0
\endpmatrix.
$$
\medbreak

\subhead
4.3.5. The cases $D^*$, $F$, $F^*$, $G$, $I$, $I^*$
\endsubhead
The cases corresponding to the nets $D^*$, $I^*$ and $F$ cannot occur.
Indeed in these cases our Relations (4.1) become
$a_1a_3,a_2a_3,a_3^2+2pa_1a_2+q a_1^2\in\M^3$ where $(p,q)=(1,0)(0,0),(0,1)$ in cases $D^*$, $E^*$ and $F$ respectively. In particular $\M^2=(a_1^2,a_2^2,a_1a_2)$. Thus
$$
a_1^2a_3=a_1a_3^2=a_2^2a_3=a_2a_3^2=a_3^3=0.
$$
Since also
$$
a_3a_j=\alpha_{3,j}a_1^2+\alpha_{3,j}^2a_2^2+\alpha_{3,j}^3a_3^2+\gamma_{3,j},\qquad j\ge4,
$$
where ${\gamma}_{3,j}\in\M^3$, $\alpha_{3,j}^h\in k$, we conclude that $a_3^2\in\Soc(A)\setminus\M^3$. 

The case corresponding to the net $F^*$ cannot occur.
Indeed in this case our Relations (4.1) become
$a_1^2,a_1a_2,a_2^2+a_3^2\in\M^3$. In particular $\M^2=(a_2^2,a_2a_3,a_1a_3)$. Thus
$$
a_1^2a_3=a_1a_2^2=a_1a_2a_3=a_1a_2^2+a_1a_3^2=0.
$$
Since also
$$
a_3a_j=\alpha_{3,j}^1a_3^2+\alpha_{3,j}^2a_2a_3+\alpha_{3,j}^3a_1a_3+\gamma_{3,j},\qquad j\ge4.
$$
where ${\gamma}_{3,j}\in\M^3$, $\alpha_{3,j}^h\in k$, we conclude that $a_1a_3\in\Soc(A)\setminus\M^3$. 

The case corresponding to the net $G$ cannot occur.
Indeed in this case Relations (4.1) become
$a_1^2,a_2^2,a_2a_3\in\M^3$. In particular $\M^2=(a_3^2,a_1a_3,a_1a_2)$. Thus
$$
a_1^2a_2=a_1a_2^2=a_1a_2a_3=a_2a_3^2=0.
$$
Since we have 
$$
a_1a_j=\alpha_{1,j}^1a_3^2+\alpha_{1,j}^2a_1a_3+\alpha_{1,j}^3a_1a_2+\gamma_{1,j},\qquad j\ge4.
$$
where ${\gamma}_{1,j}\in\M^3$, $\alpha_{1,j}^h\in k$, we infer $a_1a_2\in\Soc(A)\setminus\M^3$. 

The case corresponding to the net $I$ cannot occur.
Indeed in this case Relations (4.1) become
$a_1^2,a_1a_2,a_2^2\in\M^3$. In particular $\M^2=(a_3^2,a_1a_3,a_2a_3)$. Thus
$$
a_1^3=a_1^2a_2=a_1^2a_3=a_1a_2^2=a_1a_2a_3=a_2^3=a_2^2a_3=0.
$$
Changing possibly $a_3$ in $a_3+ua_2+va_1$ we can always assume that $\M^3=(a_3^3)$. Let $a_ha_3^2=\lambda_ha_3^3$, $h=1,2$. We have 
$$
a_3a_j=\alpha_{3,j}^1a_1a_3+\alpha_{3,j}^2a_2a_3+\alpha_{3,j}^3a_3^2+\gamma_{1,j}a_3^3,\qquad j\ge4.
$$
where ${\gamma}_{3,j},\alpha_{3,j}^h\in k$. Via $a_4\mapsto a_4+\alpha_{3,j}^1a_1+\alpha_{3,j}^2a_2+\alpha_{3,j}^3a_3+\beta_{1,j}a_3^2$, we can also assume $a_3a_4=0$. Thus if $(\lambda_1,\lambda_2)\ne(0,0)$, then $\lambda_2a_2a_3-\lambda_1a_1a_3\in\Soc(A)\setminus\M^3$. If $\lambda_1=\lambda_2=0$, then $a_1a_3\in\Soc(A)\setminus\M^3$.
\medbreak

We conclude the section with the following

\proclaim{Proposition 4.4}
The algebras $A_{n,3n+3}^{h,\alpha^2}$ are pairwise non--isomorphic for $h=1,\dots,6$.
\endproclaim
\demo{Proof}
For each $h$ we define $I_h^2$ as the ideal generated by
$$
\{\ u\in\gr(A_{n,3,n+3}^{h,\alpha^2})\ \vert\ u^2=0\ \}\cap\M/\M^2.
$$
Using our representation for $A_{n,3,n+3}^{h,\alpha^2}$, one checks that $I_h^2=(x_4,\dots,x_n)\subseteq\gr(A_{n,3,n+3}^{h,\alpha^2})$, thus $\gr(A_{n,3,n+3}^{h,\alpha^2})/I_h^2\cong\gr(A_{3^2,8}^{h,\alpha^2})$

An isomorphism $\psi\colon A_{n,3,n+3}^{h',\alpha'{}^2}\to A_{n,3,n+3}^{h'',\alpha''{}^2}$ yelds a graded isomorphism $\Psi\colon \gr(A_{n,3,n+3}^{h',\alpha'{}^2})\to\gr(A_{n,3,n+3}^{h'',\alpha''{}^2})$. We have, trivially, $\Psi(I_{h'}^2)=I_{h''}^2$, hence an isomorphism $\gr(A_{3^2,8}^{h',\alpha'{}^2})\to\gr(A_{3^2,8}^{h'',\alpha''{}^2})$. Since $A_{3^2,8}^{h,\alpha^2}$ is actually graded, it follows $A_{3^2,8}^{h,\alpha^2}\cong\gr(A_{3^2,8}^{h,\alpha^2})$, hence we infer the existence of a graded isomorphism $A_{3^2,8}^{h',\alpha'{}^2}\to A_{3^2,8}^{h'',\alpha''{}^2}$ induced by the original isomorphism $\psi$. We thus conclude that it suffices to prove the statement in the particular case $n=3$.

Notice that $A_{3^2,8}^{h,\alpha^2}$, $h\ge4$, is not a complete intersection, thus it cannot be isomorphic to $A_{3^2,8}^{h,\alpha^2}$, $h\le3$. Let $A_{3^2,8}^{h',\alpha'{}^2}\cong k[x_1,x_2,x_3]/I'$ and $A_{3^2,8}^{h'',\alpha''{}^2}\cong k[x_1,x_2,x_3]/I'$. Then $\psi$ finally induces an automorphism $\psi_0$ of $k[x_1,x_2,x_3]$ such that $\psi_0(I')=I''$, thus the corresponding nets of conics, which are generated by the generators of degree $2$ of the ideals $I'$ and $I''$ , must be projectively equivalent, hence they have isomorphic discriminant curves.

First we examine the more delicate case $A_{3^2,8}^{1,\alpha^2}$, $\alpha\in k$. We already noticed above that the discriminant curves of the corresponding nets must be isomorphic, hence they must have the same $j$--invariant or, if singular, the same kind of singular point.

In our case the net is generated by 
$$
x_1x_2+x_3^2=0,\qquad
x_1x_3=0,\qquad
x_2^2-4p x_3^2+x_1^2+2p x_1x_2=0,
$$
(here we modified the last generator of the net in order to obtain a cubic in Weierstrass form as discriminant of the net, setting $\alpha=6p$, $p\in k$)
thus its discriminant curve $\Delta$ has equation
$$
\lambda_1^2\lambda_2=(\lambda_0^2+4p\lambda_0\lambda_2+4(p^2-1)\lambda_2^2)(4p\lambda_2-\lambda_0),
$$
which is singular if and only if $\alpha=\pm1/3$ and, in this case, it carries a node.

In all the remaining cases its $j$--invariant is 
$$
j(\Delta)=-{{27p^2(1-p^2)^2}\over{(1-9p^2)^2}}.
$$
We recall that, once we fix the discriminant curve $\Delta$ of the net, there are exactly three nets of conics with discriminant curve $\Delta$ (e.g. see [Be], Chapter VI). They correspond to the three non--trivial theta--characteristics on $\Delta$. In our case there are exactly six possible values of $p$ corresponding to the same $j$--invariant for $\Delta$ and we have already checked that we can identify the two cases $\pm p$ (see Section 3.2.2). Thus we have exactly three possible values of $p$ giving rise to possible non--isomorphic nets of conics for a fixed  $j$--invariant. We conclude that such values actually correspond to non isomorphic nets of conics, thus to non--isomorphic $k$--algebras.

In a similar way, looking at the discriminant of the net (which are listed in Table 1 of [Wa]), one checks that $A_{3^2,8}^{h,\alpha^2}$ are pairwise non--isomorphic for  either $h\le3$ or $h\ge4$.
\qed
\enddemo

\remark{Remark 4.5}
We look for the deformations of the $k$--algebras $A_{n,3,n+3}^{h,\alpha^2}$. We define the following ideals in $k[b,x_1,\dots,x_n]$:
$$
\gather
J^{1,\alpha^2}:=(x_1x_2+x_3^2,x_1x_3,x_2^2-\alpha x_3^2+x_1^{2},x_ix_j,x_h^2-x_1^3,x_n^2-bx_n-x_1^3)_{{1\le i<j\le n,\ 4\le j}\atop 4\le h\le n-1},\\
J^{2,0}:=(x_1^2,x_2^2,x_3^2+2x_1x_2,x_ix_j,x_h^2-x_1x_2x_3,x_n^2-bx_n-x_1x_2x_3)_{{1\le i<j\le n,\ 4\le j}\atop 4\le h\le n-1},\\
J^{3,0}:=(x_1^2,x_2^2,x_3^2,x_ix_j,x_h^2-x_1x_2x_3,x_n^2-bx_n-x_1x_2x_3)_{{1\le i<j\le n,\ 4\le j}\atop 4\le h\le n-1},\\
J^{4,0}:=(x_2^3-x_1^3,x_3^3-x_1^3,x_ix_j,x_h^2-x_1^3,x_n^2-bx_n-x_1^3)_{1\le i<j\le n\atop 4\le h\le n},\\
J^{5,0}:=(x_1^2,x_1x_2,x_2x_3,x_2^3-x_3^3,x_1x_3^2-x_3^3,x_ix_j,x_h^2-x_3^3,x_n^2-bx_n-x_3^3)_{{1\le i<j\le n,\ 4\le j}\atop 4\le h\le n-1},\\
J^{6,0}:=(x_1^2,x_1x_2,2x_1x_3+x_2^2,x_3^3,x_2x_3^2,x_ix_j,x_h^2-x_1x_3^2,x_n^2-bx_n-x_1x_3^2)_{{1\le i<j\le n,\ 4\le j}\atop 4\le h\le n-1}.
\endgather
$$
In $k[b,x_1,\dots,x_n]$ we have
$$
J^{h,\alpha^2}=(x_1,\dots,x_{n-1},x_n+b)\cap (J^{h,\alpha^2}+(x_n^2)).
$$
Notice that in $J^{h,\alpha^2}+(x_n^2)$ there always is  a polynomial which is
$$
bx_n+\text{cubic polynomial in $x_1,\dots,x_{n-1}$}.
$$
With this in mind it is easy to check, case by case, that for a fixed $b\ne0$ in $k$,  we have $k[x_1,\dots,x_{n-1}]/(J^{h,\alpha^2}+(x_n^2))\cong A_{n-1,3,n+2}^{h,\alpha^2}$. We thus conclude that the family ${\Cal A}^{h,\alpha^2}:=k[b,x_1,\dots,x_n]/J^{h,\alpha^2}_n\to \a1$ is flat, it has special fibre over $b=0$ isomorphic to  $A_{n,3,n+3}^{h,\alpha^2}$ and general fibre isomorphic to $A_{n-1,3,n+2}^{h,\alpha^2}\oplus A_{0,1}$.
\endremark
\medbreak

\head
5. The locus $\Hilb_{d}^{G}(\p{N})$ for $d\le9$
\endhead
Taking into account the results of the previous sections, it is now possible to study the irreducibility of $\Hilb_{d}^{G}(\p{N})$ and its singular locus when $d\le9$. The first result is the following

\proclaim{Proposition 5.1}
Let $\char(k)\ne2,3$. If $d\le 9$, then $\Hilb_{d}^{G}(\p{N})=\Hilb_{d}^{G,gen}(\p{N})$ hence it is irreducible.
\endproclaim
\demo{Proof}
Let $X\in \Hilb_{d}^{G}(\p{N})$ be AS then $X\in \Hilb_{d}^{G,gen}(\p{N})$ by Proposition 2.3. We can complete the proof of the above statement if we examine the case of non--AS irreducible schemes.

It suffices to prove that each such $X\in \Hilb_{d}^{G}(\p{N})$ is a specialization of a flat family of schemes in $\Hilb_{d}^{G,gen}(\p{N})$. But these schemes are of the form $X\cong \spec(A)$ where $A$  is either $A_{n,\delta}$, with $n=4,5,6,7$ and $6\le n+2\le \delta\le 9$ or $A_{n,2,\delta}^h$ with $h=1,2$ and $8\le n+4\le \delta\le 9$ or $A_{4,3,9}^{h,\alpha^2}$ with $h=1,\dots,6$.

For example, $A_{n,\delta}$ is in the flat family ${\Cal A}_{n,d}$ (see Remark 2.9). Its general member is $ A_{n,\delta-1}\oplus A_{0,1}$ (when $\char(k)>2$) which is in $\Hilb_{d}^{G,gen}(\p{N})$ due to the argument above. The same argument for ${A}_{n,2,\delta}^h$ with families ${\Cal A}_{n,2,\delta}^h$, $h=1,2$, defined in Remark 3.4 (when $\char(k)>4\ge\lev({A}_{n,2,\delta}^h)=d-n-1$) and for ${A}_{4,3,9}^h$ with families ${\Cal A}_{4,3,9}^{h,*}$, $h=1,\dots,6$, defined in Remark 4.5 (when $\char(k)>3$), completes the proof.
\qed
\enddemo

\remark{Remark 5.2}
When $d\le7$ the above result is classically known. Indeed, in this case, $\Hilb_{d}(\p{N})$ is irreducible (see e.g. [Ma2]: see also [C--E--V--V]), hence the same is true for the open dense subset $\Hilb_{d}^{G}(\p{N})$. 

It is proved in [C--E--V--V] that $\Hilb_{8}(\p{N})$ is again irreducible if $N\le3$ and it consists of two distinct components if $N\ge4$. In this case, beside the component $\Hilb_{8}^{gen}(\p{N})$, containing all the points representing smooth schemes, there is another component. Its points represent irreducible schemes $X=\spec(A)$, where $A$ is a local and Artinian $k$--algebra with $H(A)=(1,4,3)$, thus $X\not\in\Hilb_{8}^G(\p{N})$ which again turns out to be irreducible.
\endremark
\medbreak

Now we examine $\Sing(\Hilb_{d}^{G}(\p{N}))$. To this purpose we will make use of Equality (2.2) and the comments after it. 

Let $X=\bigcup_{i=1}^pX_i$ where $X_i$ is irreducible of degree $\delta_i$. We already checked that if $X$ is AS, then it is unobstructed by Proposition 2.3.

Now assume that $X$ contains a component $Y\cong\spec(A)$ where $A$ is either $A_{n,\delta}$, with $n=4,5,6,7$ and $6\le n+2\le \delta\le 9$ or $A_{n,2,\delta}^h$ with $h=1,2$ and $8\le n+4\le9$. 

In the first case $Y$ can be deformed in $\Hilb_\delta(\p n)$ to $\widehat{Y}:=\spec(A_{n,n+2}\oplus A_{0,1}^{\oplus \delta-n-2})$ (when $\char(k)>2$). Due to Theorem 3.5 of [C--N], we have
$$
h^0\big(Y,{\Cal N}_{Y\vert \a n}\big)\ge{{(n+2)^3-7(n+2)}\over6}+n(\delta-n-2)\tag5.3
$$
hence $h^0\big(Y,{\Cal N}_{Y\vert \a n}\big)>\delta n$ for each $n\ge4$.  The same argument for ${A}_{n,2,\delta}^h$ with the families  ${\Cal A}_{n,2,\delta}^h$, $h=1,2$, yields the obstructedness of $Y$ in the second case.

The above discussion proves the following

\proclaim{Proposition 5.4}
Let $\char(k)\ne2,3$. If $d\le 8$, then $X\in\Hilb_{d}^{G}(\p{N})$ is obstructed if and only if it represents a non--AS scheme.
\qed
\endproclaim

The picture in the case $d=9$. In this case we have to examine the case of schemes $X:=\spec(A)$ where $A\cong A_{4,3,9}^{h,\alpha^2}$ with $h=1,\dots,6$. Such algebras are flat specialization of $A_{3^2,8}^{h,*}$ with $h=1,\dots,6$ which are unobstructed being AS. Thus we cannot prove (or disprove) that they are obstructed following the above method.

To this purpose we will compute directly $h^0\big(X,{\Cal N}_X\big)$ using Equation (2.2) applied to a suitable embedding $X\subseteq\p {7}$. We recall the following (see [Sch], Lemma (4.2))

\proclaim{Proposition 5.5}
Let $X$ be a Gorenstein scheme of dimension $0$ and degree $d$. Then there exists a non--degenerate embedding $i\colon X\hookrightarrow\p{d-2}$ as aG subscheme. In particular we have a resolution
$$
0\longrightarrow S(-d)\longrightarrow
S(-d+2)^{\oplus\beta_{d-3}}\longrightarrow\dots\longrightarrow
S(-2)^{\oplus\beta_1}\longrightarrow S\longrightarrow S_X\longrightarrow0
\tag 5.5.1
$$
where $S:=k[x_0,\dots,x_{d-2}]$ and
$$
\beta_h:={{h(d-2-h)}\over{d-1}}{d\choose{h+1}},\qquad h=1,\dots,d-3.
$$

Moreover if $j\colon X\hookrightarrow\p{d-2}$ is another embedding whose image is non--degenerate and aG then there exists $\varphi\in\PGL_{d-1}$ such that $\varphi\circ i=j$.
\qed
\endproclaim

We are now able to deal with the singular locus of $\Hilb_{9}^{G}(\p{N})$. Obstructed schemes $X\in \Hilb_{9}^{G}(\p{N})=\Hilb_{9}^{G,gen}(\p{N})$ are necessarily non--AS due to Proposition 2.3, thus they must contain a component isomorphic to either $\spec(A_{n,\delta})$ or $\spec(A_{n,2,\delta}^h)$, or $\spec(A_{4,3,9}^{h,\alpha^2})$.

As in the proof of Proposition 5.2 it is easy to check that in the first two cases we have obstructed schemes as soon as $n\ge4$. We examine now the third case, i.e. $X\cong\spec(A_{4,3,9}^{h,\alpha^2})$. Since the obstructedness of $X$ does not depend on the embedding, we can consider the non--degenerate aG embedding $i\colon X\hookrightarrow\p7$ of Proposition 5.5. Thus we can view at $X$ as point in $\Hilb_{9}^{G,gen}(\p{7})$ which has dimension $63$.

Splitting the sheafification of Sequence (5.5.1) into short exact sequences and taking their cohomology, we have $H^0(X,\Ofa X(t))\cong {}_tS_X$  for $ t\ge 2$, where ${}_tS_X$ is the component of degree $t$ of $S_X$. Thus, dualizing the sheafification of Sequence (5.4.1) tensorized by $\Ofa X$ and taking the associated cohomology exact sequence, we finally obtain
$$
0\longrightarrow H^0\big(X,{\Cal N}_{X\vert\p7}\big)\longrightarrow {}_2S_X^{\oplus\beta_1}\longrightarrow {}_{3}S_X^{\oplus\beta_2}
$$
which allows us to compute $h^0\big(X,{\Cal N}_X\big)$ very easily, once that the matrix of $S(-3)^{\oplus\beta_2}\to S(-2)^{\oplus\beta_1}$ is known. With this idea, using the computer software Macaulay (see [B--S]), we obtain that $h^0\big(X,{\Cal N}_{X\vert\p7}\big)=68$ if $h=4,5,6$ and $h^0\big(X,{\Cal N}_{X\vert\p7}\big)=63$ if either $h=1,2$ or $h=1$ and $\alpha=2$ (which is the unique value of $\alpha$ such that the net generated by the three conics $x_1x_2+x_3^2$, $x_1x_3$, $x_2^2-\alpha x_3^2+x_1^{2}$ has singular discriminant curve).

Since the set of $\alpha\in k$ such that $\spec(A_{4,3,9}^{1,\alpha^2})\in \Hilb_{9}^{G}(\p{N})$ is obstructed is close, we conclude we the following partial result

\proclaim{Proposition 5.6}
Let $\char(k)\ne2,3$ and $X\in \Hilb_{9}^{G}(\p{N})$. If $X\cong\spec(A_{4,3,9}^{1,\alpha^2})$, then the set of $\alpha\in k$ such that $X$ is obstructed is finite. If $X\not\cong\spec(A_{4,3,9}^{1,\alpha^2})$, then it is obstructed if and only if it contains an irreducible component isomorphic to either $\spec(A_{n,\delta})$ or $\spec(A_{n,2,\delta}^h)$, where $n\ge4$, or $\spec(A_{4,3,9}^{h,\alpha^2})$, where $h=4,5,6$. 
\qed
\endproclaim

\remark{Remark 5.7}
The proof of Proposition 5.4 can be easily generalized to prove that the closure in $\Hilb_d^G(\p n)$ of the locus $H_{n,d}$ of schemes isomorphic to $\spec(A_{n,n+2}\oplus A_{0,1}^{\oplus d-n-2})$ is always contained in $\Sing(\Hilb_{d}^{G}(\p{N}))$ for $n\ge4$. 

On one hand, this means also that the schemes $\spec(A_{4,3,9}^{h,\alpha^2})$, $h\le3$ are not in $H_{4,9}$. On the other hand, let us take $X\in H_{4,9}$ and embed it in $\p7$ as non--degenerate aG subscheme. Then Equations (2.2) and (5.3) easily yield $h^0\big(X,{\Cal N}_{X\vert\p7}\big)\ge68$, thus we cannot exclude that $\spec(A_{4,3,9}^{h,\alpha^2})\in H_{4,9}$, $h=4,5,6$.

Indeed this is actually the case for $h=4$. To prove this consider the ideal $J\subseteq k[b,x_1,\dots,x_n]$
$$
\align
J:&=(x_2^3-bx_1^2-x_1^3,x_1^3-bx_3^2-x_1^3,x_ix_j,x_h^2-bx_1^2-x_1^1)_{1\le i<j\le n\atop 4\le h\le n}=\\
&=(x_1+b,x_2,\dots,x_n)\cap(x_2^3-bx_1^2,x_1^3-bx_3^2,x_ix_j,x_h^2-bx_1^2)_{1\le i<j\le n\atop 4\le h\le n}.
\endalign
$$
Thus ${\Cal A}:=k[b,x_1,\dots,x_n]/J\to \a1$ is a flat family having special fibre over $b=0$ isomorphic to  $A_{n,3,n+3}^{4,0}$ and general fibre isomorphic to $A_{n,2,n+2}^{1}\oplus A_{0,1}$. Thus, for each $n\ge3$, we have $X:=\spec(A_{n,3,n+3}^{4,0})\in H_{4,n+3}$.

It is then natural to state the following

\definition{Conjecture 5.7.1}
$\Sing(\Hilb_{d}^{G,gen}(\p{N}))=H_{4,d}$ for each $d$.
\enddefinition
\endremark
\medbreak

\enddocument